\input amstex
\input xy
\xyoption{all}
\documentstyle{amsppt}
\document
\magnification=1200
\NoBlackBoxes
\nologo
\hoffset1.5cm
\voffset2cm
\pageheight {16cm}

\def\C{\bold{C}}

\def\K{\bold{K}}
\def\k{\bold{k}}
\def\P{\bold{P}}
\def\M{\Cal{M}}
\def\O{\Cal{O}}
\def\T{\Cal{T}}
\def\U{\Cal{U}}



\centerline{\bf MIRRORS, FUNCTORIALITY, AND DERIVED GEOMETRY }

\bigskip

\centerline{\bf Yuri I.~Manin}

\medskip

\centerline{\it Max--Planck--Institut f\"ur Mathematik, Bonn, Germany}

\bigskip

\centerline{\bf Contents}

\medskip

1. Introduction: a mystery of Quantum Cohomology

\smallskip

2. Deformation functors and controlling DGLAs

\smallskip

 3. Deformations of analytic local rings and mirror phenomena
 
 \smallskip
 
4. Extended deformation functors and controlling $L_{\infty}$--algebras

\smallskip

5.  $F_{\infty}$--structures on extended deformation spaces

\smallskip

\quad  References

\medskip

{\it Abstract.} In this survey, I suggest to approach the problem of functorial
properties of quantum cohomology by drawing lessons from several versions
of Mirror duality involving deformation spaces.

\bigskip

\centerline{\bf 1. Introduction: a mystery of Quantum Cohomology}

\medskip

{\bf 1.1. A brief summary.}  Moduli spaces/stacks $M$ of stable curves of all genera with a finite number of marked 
points
endowed with natural correspondences between them 
form a (modular) operad $h(M)$: see [KoMan],  [BehMan], [Man1],
and subsequent works.

\smallskip

This operad acts upon
each smooth complete algebraic variety/DM--stack $V$ via correspondences.
Thus, in a wide sense of the word, motive/cohomology $h(V)$ of $V$
is endowed with 
a canonical structure of {\it algebra over $h(M)$. }
This structure is called  {\it Quantum Cohomology (QC)} of $V$.
\smallskip

{\bf A  mystery}: unlike motives/cohomology theories,
we practically do not understand 
properties  of QC considered as {\it functor of} $h(V)$.
\smallskip
A related mystery: self--referentiality of the operad $h(M)$, i.~e. its
action
upon its own components, and interaction of it with operadic structure. This problem was explicitly addressed in [ManSm1], [ManSm2].

\smallskip
Here I suggest to approach this problem using 
certain constructions
traditionally used in one of the contexts of Mirror Symmetry: {\it deformation spaces}
and their {\it enriched/derived versions.}

\smallskip

The point is that these constructions are intrinsically and richly functorial, so that one approach
to our mystery consists in bringing them back from the Looking Glass Land. 

\medskip
{\bf 1.2. Recollection and notation.} Here are somewhat more precise notations and statements. Quantum Cohomology of an irreducible smooth projective manifold $V$
is  the system of motivic morphisms
$$
I^V_{g,n,\beta}:\, h(V^{n}) \to h(\overline{M}_{g,n}).
$$

Here $\overline{M}_{g,n}$ denotes the moduli DM--stack of stable curves of genus $g$ with
$n$ marked 
points, $h$ denotes the respective motive, and $\beta$ runs over divisor classes of $V$. 
\smallskip
 This system expresses the canonical action of the motivic modular operad
 upon every ``total'' motive $h(V):$ I use here the word "total" in order to
 stress that it is not clear at all upon which direct summands of total motives
 this operad acts.

\smallskip

In the framework of this survey, I am focusing on the case $g=0, n\ge 3$,
and consider only that part of information about this action which is 
compressed in the
{\it genus zero  quantum cohomology ring $H^*_q(V)$. } 

\smallskip

Assume for simplicity that $V$ is defined 
over a field of characteristic $0$, and denote by $H^*(V):=H^*(V,K)$ its cohomology 
ring with coefficients in a $\bold{Q}$--algebra $K$.

\smallskip
As a graded $K$--module, this ring is  free of finite rank; let  $(\Delta_a), a=0,1,\dots ,r$, 
be its free
graded basis such that $\Delta_0$ is the (dual) fundamental class 
of $V$, that is,
the identity of the local artinian ring $H^*(V)$.

\smallskip
The dual homology module $H_*(V,K)$
can be considered then as the module of 
linear coordinates  on $H^*(V,K)$
with graded coordinates $(x_a)$, dual to  $(\Delta_a)$.   
We replace those $x_a$ for which $\Delta_a\in H^2(V,K)$ by their formal
exponents $q_a=e^{x_a}$ 
and construct the ring of formal series
$K_q:=K[[q_a; x_b]]$ (Novikov's ring) 
where $b$ runs over subscripts for which $\Delta_b\notin H^2(V,K)$.

\smallskip
{\it The genus zero quantum cohomology ring}  $H_q^*(V)$ is then the free $K_q$--module
$K_q\otimes_K H^*(V)$ with graded commutative multiplication which is the deformation
of the multiplication in $H^*(V)$ in the following sense:
\medskip
(i)  $\Delta_0$ remains the identity
in the deformed ring. 
\medskip
(ii) Modulo the maximal ideal $(q^a, x_b\,|\,b\ne 0)$, the ring structure of $H_q^*(V)$ is 
the same as that of $H^*(V)$. In other words, the ring $H^*_q(V)$ is a formal deformation
of the ring $H^*(V)$.

\medskip

(iii) Finally, the deformed (``quantum'') multiplication $\circ$ itself has the following 
structure. 

\smallskip
One starts with constructing
the {\it potential} $\Phi \in K_q$ whose coefficients are genus zero 
Gromov--Witten invariants expressing
(appropriately defined virtual) numbers of rational 
curves in $V$ with marked points restricted by incidence conditions.
Then one constructs the third derivatives $\Phi_{abc}:=\partial_a \partial_b \partial_c\Phi$
where $\partial_a  :=\partial/\partial x_a$. 
And finally one sets
$$
\Delta_a\circ \Delta_b := \sum_{cd}\Phi_{abc}g^{cd}\Delta_d
\eqno(1.1)
$$
where $(g^{ab})$ is the matrix of the Poincar\'e duality product upon $H^*(V)$.
For more details and examples, cf.~ [BaMan].

\medskip

{\bf 1.3. Mirrors and functoriality.} In the main body of the paper,  this concrete 
{\it quantum cohomology deformation}
of $H^*(V)$ will be considered in the more general  context of various
{\it deformation\ theories.}   These deformation theories, on the one hand,
are used in several of the many Mirror Symmetry constructions, and on the other
hand, 
these theories, especially
their extended and derived versions, have very rich
functorial properties.  

\smallskip

Functoriality of deformation spaces is explained in section 2 of this survey:
it is achieved by putting deformation theories into the general
framework of {\it controlling DGLAs} (Differential Graded Lie Algebras).

\smallskip

The section 3 is dedicated to concrete examples of applications
of deformation theories in Mirror Symmetry constructions.
Finally, sections 4 and 5 considerably extend the framework
of controlling DGLAs by introducing derived geometric constructions.

\smallskip

In order to trace back  historical roots and subsequent
developments of deformation philosophy, I can suggest the sources [Gr], [Ar], [Art1], [Art2], [Dr],
[MatY], [Sc], 
[Mi], [BuMi], and [Ma1]--[Ma3].

\smallskip

Here is a brief introduction into
one of the most basic cases and its treatment as a Mirror Symmetry construction.

\medskip

{\bf 1.4. Deformations of local artinian rings:
the case of Jacobi rings of isolated hypersurface
singularities.}  
Let $f=f(x_1,\dots ,x_m)$ be the germ of  holomorphic (or formal) function $f:\, (\C^m,0)\to (\C,0).$
Assume that $(0)$ is the only critical point of this germ.
Let 
$J(f):= \O_{\C^m,0}/(\partial f /\partial x_k)$ be its Jacobi algebra. The number $\mu$, its  linear dimension
over $\C$,  is called {\it the Milnor number} of the singularity.

\smallskip

An {\it unfolding} (or deformation) of $f$ is a holomorphic germ $F(x_1,\dots ,x_m; t_1,\dots ,t_n)$ 
at $(\C^{m+n},0)$ such that 
$F(x;0)=f(x)$. Its {\it base} is the germ $M=(\C^n, 0)$, with coordinates $(t_j)$.
We will say that a germ of tangent vector field to $(\C^{m+n},0)$ is {\it vertical}
if its projection to the base $(\C^{n},0)$ vanishes.
\smallskip
Such an unfolding is called {\it versal}, if any other unfolding
can be induced from it 
by an appropriate morphism of bases, and {\it semiuniversal}
if it is versal and 
its base has {\it the minimal dimension.} This is the first explicit expression of functoriality
that we met in this survey.

\smallskip
For more details about morphisms of unfoldings involved here, see [He], pp. 62--63.
A version of this definition is discussed on page 64 of [He].

\smallskip

Here is a criterium for checking (semiuni)versality of an unfolding $F$
(cf. [He], Theorem 5.1, and references therein). 
\smallskip
Consider the {\it critical space} $C=C_F$
of the map $(F, pr_M):\, (\C^m\times M,0)\to (\C\times M,0)$ 
which is defined
by the ideal $J(F):= (\partial F/ \partial x_k)$. For a germ of tangent vector field 
$X\in \T_{M,0}$,
denote by $\widetilde{X}$ a lift of $X$ to $(\C^m\times M, 0).$ 

\smallskip
Since the difference of any two lifts of the same $X$
must be vertical, its restriction upon $C$ vanishes 
so that the map
$$
\T_{M,0}\to  pr_{M*}\O_C :\quad   X\mapsto \widetilde{X}F\,\roman{mod}\,J(F)
$$
is well defined.
\medskip
{\bf 1.4.1. Theorem} ([He], p.~63).
{\it a) An unfolding $F$ of an isolated hypersurface singularity is 
versal iff the map 
$X\mapsto \widetilde{X}F\,\roman{mod}\,J(F)$  is surjective.

\smallskip

b) An unfolding $F$ of an isolated hypersurface singularity is 
semiuniversal 
iff the map $X\mapsto \widetilde{X}F\,\roman{mod}\,J(F)$ is bijective.}

\medskip
As soon as a semiuniversal unfolding is chosen, we can define a commutative
associative $\O_M$--bilinear multiplication $\circ$ on $\T_{M,0}$ by simply lifting it from
$pr_{M*}\O_C$ i.~e., by putting
$$
\widetilde{X\circ Y}F \,\roman{mod}\,J(F) :=\widetilde{X}F\cdot \widetilde{Y}F\,
\roman{mod}\,J(F).
\eqno(1.2)
$$

{\bf 1.5. Example: semiuniversal unfolding
of the singularity} $A_r$ {\bf vs. quantum cohomology ring} $H_q^*(\P^r).$
One of the versions 
of mirror symmetry 
 starts with the observation that  if 
a classical cohomology 
ring $H^*(V)$ is isomorphic to the Jacobi ring $J(f)$ of an isolated
hypersurface singularity, 
then locally (or formally) near the initial point of semiuniuversal
unfolding space 
the formal spectrum of
 the ring $K_q=K[[q_a;x_b]]$ 
defined in sec.~1.2 must admit 
a natural (generally non--unique) map to the germ of this unfolding,
such that 
the relative formal spectrum  of $H_q^*(V)$ is induced by the relative spectrum of
the critical  space $C_f$ defined above.

\smallskip

Moreover, one should then try to constrain a choice of this map by requiring 
its compatibility with additional 
structures induced on the unfolding spaces of two mirror sides.
The first and most important is the compatibility of the Frobenius--multiplication (1.1)
with $F$--multiplication (1.2). Besides, one should try to transfer the  
canonical grading and flat structure
on the quantum cohomology side to the unfolding side where they are 
initially absent. 

\smallskip
All these details were thoroughly studied for many homogeneous spaces of classical Lie/algebraic groups:
see a recent report [GoPe], anf for physical motivation [BelGeKo]. 
The simplest example is a  projective space $\P^r$ over $\C$.
\smallskip
Its classical cohomology with coefficients in any ring $K$  is the free $K$--module
freely generated by  $\Delta_a$ where $\Delta_a$ is the dual class of
$\P^{r-a}\subset \P^r.$ As a ring, 
$H^*(\P^r, K)$ is thus canonically isomorphic to
$K[\Delta_1]/(\Delta^{r+1}_1)$, with $\Delta_a \cong \Delta_1^a \, \roman{mod}\, \Delta_1^{r+1}$.

\smallskip

On the other side, the germ of the function of one variable $f(x)=x^{r+2}$ 
has the same Jacobi ring 
$K[x]/(x^{r+1})$. Its semiuniversal unfolding space is  
the affine space of coefficients of
the polynomial 
$$
p(t):=x^{r+1}+ t_1x^{r-1} + \dots + t_{r-1}x +t_r
$$
For many further details cf.~ [Man1], Ch.~II, sec.~4; Ch.~I, sec.~4.
In particular, it is explained, how to introduce additional flat structure and grading
that make it compatible with the respective structures upon quantum cohomology
of $\P^r$. 

\smallskip

{\bf 1.6. Homogeneous spaces.}  This approach to the quantum cohomology of $\P^r$ and its
mirror was (at least partially) extended to more general homogeneous spaces $G/P$.
The question for what such spaces  their classical cohomology is isomorphic to the
Jacobi ring of an isolated singularity, seemingly does not have a direct answer
in the literature. Probably, the answer is positive at least for minuscule/cominiscule
homogeneous spaces. In any case, the ample known information about 
explicit descriptions of cohomology rings was used in order
to successfully produce also a description of their quantum cohomology
focused more on the flat structure and structure connections than upon
$F$--multiplication, cf. [Sa].

\smallskip

For some recent surveys/original results see [ChMPe], [LamTem],
and references therein. See also an interesting extension of 
this method in [GoPe], where the authors construct deformations of Jacobi
rings of polynomials, and in the context of homogeneous spaces,
apply these constructions to equivariant cohomology and $K$--theory.

\smallskip

{\bf 1.7. Extended deformations and derived geometry.} The sections 4 and 5 of this survey will be dedicated to the problem of extending
deformation contexts, if the basic theory involving controlling DGLAs
is not satisfying enough for discussing mirror phenomena.
In particular,

\smallskip

{\it  What to do when $H^*(V)$ is NOT a Jacobi ring?}

\smallskip

For starting steps, see articles  [Mi], [SchSt], where the deformation theory of $H^*(M)$ is contained, restricted to
those deformations that, for complex compact $M$,  deform only the complex structure. 
In this spirit, finite--dimensional graded Artin rings are also considered.

\smallskip

More general is the suggestion to use {\it not} the naive
deformation spaces, but some higher step of the ladder involving the so called {\it extended}
deformation spaces/functors. We present a survey of several
of  lower steps where the functoriality of main definitions and constructions 
is the primary concern. Of course, this means that derived and higher derived versions
of all objects involved should be briefly presented. In particular, we must pass

\bigskip

-- From operads classifying these objects to 
complexes/simplicial sets/... 
of operads  to homotopy 
via model structures on the respective categories.

\medskip

-- From categories to 2--categories to $\infty$--categories ...

\medskip

-- From ``affine'' objects to gluing to ...

\medskip

Much more  details are given in the monograph [LoVa].
See also [Pr], [DoShVa1], [DoShVa2], [Lu], [To]. In particular,
To\"en's survey is a magnificent introduction to the derived
deformation theories.

\newpage

\centerline{\bf 2. Deformation functors and controlling DGLAs}

\medskip

{\bf 2.1. Formal deformation philosophy.} Let $M$ be a ``space'' that
in the next few paragraphs will embody an intuitive idea of
``space of deformations of certain structured object $X$''. Thus,
$X$ itself will correspond to a point $x\in M$. 

\smallskip

In the formal deformation philosophy, we want to get hold of 
``infinitesimal neighbourhoods'' of $x$ in $M$, or even of ``germ of $M$''
at $x$ if we can speak about analytic moduli spaces.

\smallskip

Imagine first the simplest case when $M$ (or a neighbourhood of $x$ in $M$)
is a scheme defined over a field  $\K$, and  $x$ is a $\K$--point of $M$.
Then a basis of infinitesimal neighbourhoods of $x$
consists of affine spectra of rings $\O_x/m_x^n$, $n=1,2,3, \dots$,
tautologically embedded into $X$. Here $m_x$ is the maximal ideal of $\O_x$.

\smallskip

In order to construct  this basis, or its version, we must understand
deformations $X_A$ of $X$ over local Artin $\K$--algebras $A$, but now
up to an isomorphism over $A$, so we have to consider {\it groupoids
of deformations} over variable bases, forming 
 a contravariant functor $Art_{\K} \to Grpd$
which is a basic example of  deformation functors.

\smallskip

In fact, our initial view of algebraic geometry adjusted to moduli
problems can and ought to be vastly extended, including spectral
geometry, various versions of derived geometry etc.
This is explained in [Lu], but we will not try to explain it in this 
short survey, cf. [To].

\medskip

{\bf 2.2. Deformation functors.} Historically,
early abstract theories of deformation functors were developed in [Art1], [Art2],  [Gr], [Sch], [SchSt].

\smallskip

Here we start with Deligne's dense expression  {\it controlling DGLA}, that was explained
in a letter of Pierre Deligne [De]. According to a quotation from [GoMi], Deligne
observed that  
\smallskip
{\it ``in characteristic zero, a deformation problem is controlled by a differential graded Lie
algebra, with quasi--isomorphic differential graded Lie algebras giving the same
deformation theory.''} 
\smallskip
A choice of controlling DGLA provides another construction
of a functor, and identification of both versions furnishes strong
tools for studying deformations.

\smallskip

In the remaining part of this section, we focus on the construction of the (potential) deformation functor from the controlling
DGLA, mostly accepting the conventions of [GoMi],
and later return to the task of identifying two constructions.

\medskip

{\bf 2.3. Maurer--Cartan equation}. Let $\K$ be a field of characteristic
zero and $g=(\oplus_{i\ge 0} g^i, d)$  a DGLA over $\K$, with $d:\,g^i\to g^{i+1}$, $d^2=0.$
Skew--commutativity and Jacobi identities are also supposed to be graded, with Koszul's signs.

\smallskip
The set of {\it Maurer--Cartan elements} of $g$
is defined as
$$
MC(g):= \{x\in g^1\,|\, dx+\frac{1}{2}[x,x]=0\}.
\eqno(2.1)
$$
This Maurer-Cartan equation is equivalent  to the  flatness of the respective connection on $g$ 
that sends any $f\in g^0$ to $\nabla_x(f):= df+[x,f]$, namely $\nabla_x^2=0$.
\smallskip

We now want to  identify those elements of $MC(g)$ that are connected
by the flow corresponding to this action of $g^0$ (see a more sophisticated version of this identification
in the next subsection). In order to do it properly,
one can assume as in [GoMi] (p.~48) that $\K$ is $\C$ or $\bold{R}$, and consider
 the action of simply connected Lie group with Lie algebra $g^0$.
 
 \smallskip
 Another version, not involving restrictions on $\K$, assumes instead that
 $g^0$ is nilpotent, together with its action upon $g^1$, or even that $g$ is nilpotent: see [LoVa],
 p.~499. Then one can construct the respective nilpotent algebraic group
 and its action. In each of these cases, the standard formula for the action
 of one--parametric subgroups is applicable: for $a\in g^0,$ $e^{ta}$ sends $x\in MC(g)$
 to
 $$
 e^{t\, \roman{ad}\,a}(x) + \frac{\roman{id} -e^{t\, \roman{ad}\,a}}{\roman{ad}\,a} (da) .
 \eqno(2.2)
 $$
In the most important for us series of examples, we start with arbitrary DGLA $g$
and finite dimensional nilpotent commutative $\K$--algebra  $m$ (so that $\K\oplus m$
is a local Artin algebra with maximal ideal $m$). Then $g\otimes m$ is nilpotent,
with grading and $d$ coming from $g$.

 \medskip

  {\bf 2.4. Philosophy of controlling DGLAs.} The general scheme is as follows: starting with a chosen
  deformation problem we construct {\it groupoids} and
  {\it arrows} in the following diagram:
  
  \medskip
  
{\it  Deformation groupoid  $\Longrightarrow$ Controlling DGLA $L$ $\Longrightarrow$ 
 Groupoid associated to $L$}
  
  \medskip
  
  and finally establish an equivalence between two groupoids in it.

\smallskip

An explicit construction of the relevant DGLA (first arrow here) requires creative thinking and 
the study of instructive examples,
existing in the literature. The second arrow is  somewhat more
standardised, and we will start with it.

\medskip
  {\bf 2.5. Deligne groupoid $D(g,A)$. } Let $g$ be a DGLA as above and $A=\K\oplus m_A$
  an Artin local algebra.
  
  \smallskip
  
  Then we put
  $$
  \roman{Ob}\, D(g, A) := MC(g\otimes m_A)
  \eqno(2.3)
  $$
  and for $x,y \in MC(g\otimes m_A)$
  $$
  \roman{Hom} (x,y) := \{ a\in g^0\otimes m_A \,|\, e^{a} (x)= y\}.
  $$
  Finally, the composition of morphisms is defined via (2.2).
  
  \smallskip
  
  An elementary, but important remark is that $D(g,A)$  is itself
  a covariant functor of $(g, A)$ considered as a variable object
  of the categorical product of DGLAs with $Art_K$. More precisely ([GoMi], p.~53):
  
  \smallskip
  
  (i) For any homomorphism of DGLAs  $\varphi :\, g\to h$ there is a natural
  functor $\varphi_*:\, D(g,A)\to D(h,A)$.
  
  \smallskip
  
  (ii) For any homomorphism of Artin local  $\K$--algebras $\psi :\,A\to A^{\prime}$
there is a natural
  functor $\psi_*:\, D(g,A)\to D(g,A^{\prime})$.
  
  \smallskip
  
  (iii) These functors can be chosen in such a way that for $(\varphi , \psi ):\, (g,A)\to (h,A^{\prime})$
  we have  the equality (and not just an equivalence) of functors
  $\psi_*\varphi_* = \varphi_*\psi_* :\, D(g,A)\to D(h,A^{\prime}).$

  \smallskip
  
  The critically important property of this construction is this: if $\varphi$
  is a quasi--isomorphism of DGLAs, then $\varphi_*$ is an equivalence of groupoids.
  Actually, for the construction of $MC (g\otimes m_A)$ only $g^i$ with
  $i=0,1,2$ are essential, so that we have ([GoMi], Theorem 2.4):
  
  \medskip 
  
  {\bf 2.5.1. Proposition.} {\it  If $\varphi$ induces isomorphisms $H^i(g)\to H^i(h)$
  for $i=0,1$  and a monomorphism for $i=2$, then $\varphi_*$ is an equivalence
  of groupoids. }

  \bigskip
  
  \centerline{\bf 3. Deformations of analytic local rings and mirror phenomena}
  
  \medskip

{\bf 3.1. Groupoids associated to deformations of analytic local rings.} Here we will
illustrate on a concrete example  both steps involved in realisation of the philosophy
briefly sketched in sec.~3.5. above. For a detailed treatment of this example,
see [BuMi], sec.~5.

\smallskip

Let $\k$ be a complete normed field of characteristic zero.
Denote by  $\k \langle z_1, \dots ,z_m \rangle $ the ring of convergent power series in $(z_k)$.
{\it An analytic local $\k$--algebra $B$} is a quotient of $\k \langle z_1, \dots ,z_m \rangle $
modulo a (topologically closed) ideal. Denote by $Art_{\k}$ the category of
Artin local $\k$--algebras. 
\smallskip

Now fix an analytic local $\k$--algebra $B$. 
\smallskip

Below I essentially use intuition and conventions related to the version of definition of
a moduli groupoid  explained in [Man1], Ch.~V, Sec.~3.1 and 3.2, pp.~210--211.
One notational difference is that since we deal here with affine schemes and/or their versions, omitting
the passage to their (Grothendieck) spectra,  arrows in the respective categories  
are inverted in comparison with those in [Man1]. The adjective  {\it ``cofibered''} below reminds about this.

\smallskip
{\bf  3.1.1. Definition.} The cofibered groupoid $\roman{Def} (B)$ of the deformations of $B$ consists of the following data:

\medskip

{\it Category of bases.} This is the category $Art_{\k}$. For any object $A$
of this category, we denote by $m_A$ its maximal ideal.

\smallskip

{\it Category of families.} One object $(B^{\prime},\rho )$ of this category $\roman{Def} (B;A)$
(intuitively, a  family over the base which is the  spectrum of
$A\in Ob\, Art_{\k}$)  consists of a flat $A$--algebra $B^{\prime}$
and a morphism of $A$--algebras $\rho :\, B^{\prime} \to B$
which induces an isomorphism $\bar{\rho} :\, B^{\prime}/m_A B^{\prime} \to B.$

\smallskip

One morphism $(B^{\prime},\rho_1)\to (B^{\prime\prime},\rho_2)$ in $\roman{Def} (B;A)$ is a homomorphism of
$A$--algebras $\varphi :\, B^{\prime}\to B^{\prime\prime}$ which modulo $m_A$ induces
the identity morphism of $B$.

\medskip

{\it Base change functor.} Given a morphism $A_1\to A_2$, the respective base change functor
$\roman{Def} (B;A_1) \to \roman{Def} (B;A_2)$ is $*\mapsto A_2\otimes_{A_1}*$
where $*$ stands for respective objects, morphisms and diagrams.

\medskip

{\bf 3.1.2. Lemma.} {\it All endomorphisms in $\roman{Def} (B)$ are isomorphisms. Moreover, 
they are exponentials of nilpotent derivations.} ([BuMi], p. 45.)

\medskip

{\bf 3.2. Passage to resolutions of $B$.} Let $R^{\bullet}$ be a free graded commutative $\k$--algebra, with
$R^m=0$ for $m>0$,
endowed with a differential $\partial$ of degree one and a surjective homomorphism
$\varepsilon :\, R^{\bullet}\to B$ which is a quasi--isomorphism. Then $(R^{\bullet}, \partial , \varepsilon )$
is called  {\it a multiplicative resolution, or resolvent} of $B$ over $\k$.

\smallskip

Sometimes it is convenient to work instead with $R_{\bullet}$ where $R_{m}=R^{-m}$.

\smallskip

Such resolutions exist and are unique up to homotopy equivalence.

\medskip

{\bf 3.2.1. Definition.} Let $R$ be a resolution $R^{\bullet}$ as above.  The  groupoid $\roman{Def} (R)$ of deformations of $R$
cofibered over $Art_{\k}$ consists of the following data:

\medskip

 {\it Category of bases} remains to be $Art_{\k}$.
 
 \smallskip

 {\it Category of families.}  One object $(R^{\prime},\rho )$ of this category $\roman{Def} (R;A)$
(intuitively, a  family of resolutions over the base which is the  spectrum of
$A\in Ob\, Art_{\k}$)  consists again of several components. 

\smallskip

The first one is a flat deformation $R^{\prime}$ of the algebra $R$ over $A$. Since
$R$ is free, we may and will henceforth assume that $R^{\prime} = R\otimes_{\k}A$
so that $R^{\prime}/m_A R^{\prime}= R.$

\smallskip

The second component is a differential $\partial^{\prime}$ of $R^{\prime}$ deforming 
$\partial$.

\smallskip

\smallskip

{\bf Fact} ([BuMi], p. 46). $R^{\prime}$ is a resolution of $H_0(R^{\prime})$
by free $A$--modules.

\medskip

{\it One morphism} $\varphi :\,(R,\partial )\to   (R^{\prime},\partial^{\prime} )$ is a homomorphism
of differential graded algebras such that $\varphi \equiv id \, \roman{mod}\, m_A$.

\medskip

{\it Base change functor.} It is again $*\mapsto A_2\otimes_{A_1}*$
where $*$ stands for respective objects, morphisms and diagrams.

\medskip

{\bf  3.3. Controlling DGLAs.} They will belong to a general class of DGLAs
defined in [BuMi], p.~4 in the following way. 

\smallskip
Let $V=V^{\bullet}$ be a non--negatively graded vector space over $\k$.
An endomorphism $T$ of $V$ of degree $l$ is a  linear map
$T:\, V^{\bullet}\to V^{\bullet +l}$. The space of such maps is denoted $Hom^l(V,V)$.
Their {\it direct sum} is denoted $Hom\, (V,V)$. It is a {\it graded Lie (super)algebra} with
commutator
$$
[S,T]:=S\circ T- (-1)^{ij} T\circ S
$$
for $S$, resp $T$, of degree $i$, resp. $j$.

\smallskip

Now assume that $V$ is a graded commutative algebra. Denote by $Der\,V$ the space of its graded derivations
over $\k$. It is closed wrt [,]. Start with this algebra, or usually its Lie subalgebra of non--negative degree
$Der^+(V)$.

\smallskip

Usually our DGLAs will be $L=Der^+V$ {\it endowed with an additional derivation $d:\,L\to L$
of degree 1 with $d^2=0$.}

\medskip

{\bf 3.4. Groupoids associated to DGLAs.}
Let $L=(L^{\bullet} , d)$ be a DGLA. We associate with $L$ its deformation groupoid $C(L)$ cofibered over 
$Art_{\k}$.

\medskip

{\bf 3.4.1. Definition.} The groupoid $C(L)$ consists of the following data:

\medskip

{\it Category of bases.} It is  $Art_{\k}$.

\smallskip

{\it Category of families over $A$: objects.} One object of the category  $C(L;A)$
is an element $\eta \in L^1\otimes_{\k} m_A$ satisfying the equation
$$
d\eta +\frac{1}{2}[\eta ,\eta ] = 0.
$$

 {\it Category of families over $A$: morphisms.} In [BuMi], p. 5, morphisms are
 defined in the following way:
 $$
 Mor\, C(L;A) := exp (L^0\otimes m_A).
 $$
 Here  $exp (L^0\otimes m_A)$ is a nilpotent Lie/algebraic group
 with underlying space $L^0\otimes m_A$ and Campbell--Hausdorff composition
 $$
 X\cdot Y := log\ (exp (X) exp (Y)).
 $$

  The morphisms act on objects by the ``affine action'':
 $\lambda \in L^0\otimes m_A$ sends $\eta \in Ob\, C(L;A) \subset L^1\otimes m_A$ 
 to $\alpha (e^{\lambda}\cdot \eta )$. The latter element is
 determined by the formula
 $$
 d\alpha (\lambda )\cdot \eta = [\lambda ,\eta ] - d\lambda .
 $$
 
 \medskip
 
 A slightly more transparent version is given in [R--N], p.~2. Each element $\lambda \in L^0\otimes m_A$
 defines a ``vector field'' on $L^1\otimes m_A$ sending $\eta \in L^1\otimes m_A$
 to
 $$
 d\lambda + [\lambda,\eta ]\in L^1\otimes m_A.
 $$
 It is tangent to the Maurer--Cartan locus in the following sense:  
 if $\eta (t)$ is a flow of $\lambda$, that is
 $$\frac{d}{dt} \eta (t) = d\lambda + [\lambda, \alpha (t)]
 $$
 with $\eta (0)$ satisfying Maurer--Cartan, the all $\eta (t)$ satisfy
 Maurer--Cartan.
 
 \smallskip
 
 Then the set of morphisms $\eta_0\to \eta_1$ is defined as the set of $\lambda\in L^0\otimes m_A$
 such that the flow starting with $\eta_0$ for $t=0$ produces $\eta_1$ for $t=1$.

\medskip

{\it Base change functor.} It is again induced by $*\mapsto A_2\otimes_{A_1}*$
where $*$ stands for respective objects, morphisms and diagrams.

\medskip

{\bf 3.5. Equivalence of deformation groupoids and DGLA groupoids}
(see  [BuMi], pp. 47--48).
Let again $B$ an analytic local $\k$--algebra, $R= (R^{\bullet}, \partial )$ its resolution as above,
$L= (L^{\bullet},d)$ be its {\it tangent complex}:  the differential graded Lie algebra of graded derivations of $R$ of non--negative degree,
and $d:= \roman{ad}\,\partial$. This means that for $\eta\in L^i$
$$
d\eta = \partial\circ \eta - (-1)^i \eta\circ \partial .
$$
\medskip

We wish now construct an equivalence of groupoids $p:\, C(L)\to \roman{Def} (R).$ Let $A$ be a local artinian $\k$--algebra.
For an object of $C(L;A)$, $\eta \in L^1\otimes m_A$, we must first of all define its image
as an object of $\roman{Def} (R;A)$. Recall that an object of $\roman{Def} (R;A)$ is represented by a flat
differential graded $A$--algebra $(R^{\prime}, \partial^{\prime})$ and a
map $\rho : R^{\prime} \to R$. In particular, we may and will assume that $R^{\prime}=R\otimes A$.

\smallskip
Recall that a morphism in $\roman{Def} (R;A)$ is a homomorphism of 
graded commutative algebras reducing  to identity modulo $m_A$.
Denote by $\beta$ the canonical isomorphism ($Hom$ overlooks the differentials)
$$
\beta:\ Hom_{\k}(R,R)\otimes A\to Hom_A(R^{\prime}, R^{\prime}),\  \beta (\eta\otimes t)= t(\eta\otimes id).
$$

\smallskip

Now we can define $p(\eta )$ for any object $\eta$ of $C(L;A)$, that is, $\eta\in L^1\otimes_{\k}m_A$
satisfying the Maurer--Cartan equation:
$$
p(\eta ):= (R\otimes A, \beta (\partial )+ \beta (\eta )).
$$

Finally, we can define $p$ on morphisms: for $\roman{exp} (\lambda ) \in  Mor C(L;A)$ we put
$$
p (\roman{exp} (\lambda )) := \beta  (\roman{exp} (\lambda )).
$$

{\bf 3.5.1. Claim} ([BuMi], p. 47). The functor $p$ is an equivalence of groupoids.

\medskip

{\bf Comparison of groupoids $h: \roman{Def} (R) \to \roman{Def} (R).$} 
This functor is defined on objects $R^{\prime}$, resp. morphisms $\varphi$, by
$$
h(R^{\prime}) := H_0(R^{\prime}),\quad h(\varphi ):= H_0(\varphi ).
$$

\smallskip

{\bf 3.5.2. Claim} ([BuMi], p. 53). $h$ induces an isomorphism of functors
$$
h: \roman{Iso\, Def (R)} \to  \roman{Iso\, Def (B)}
$$
where $\roman{Iso}$ are the sets (or small categories, p. 5) of isomorphism classes.
\medskip

{\bf 3.6. Examples: Mirror symmetry in the Looking Glass Land.} Here we briefly describe
a version of Mirror Symmetry in which {\it both sides} are Deligne (Maurer--Cartan) groupoids associated with
{\it different} DGLAs: see
[ClLaPo], pp. 4--6, [ClOvPo], and [ClPo].  In these examples, the central role is played by
an additional structure on the controlling DGLAs which is introduced from the start:
namely, they are Differential Graded Gerstenhaber Algebras, or briefly DGAs:
see [Po] for a very detailed description.

\smallskip

Let $(h,J)$ be a real Lie algebra with integrable
complex structure on it.  Starting with this datum, one can define a controlling $DGA(h,J).$ Let now $k$ be a real Lie algebra 
with a symplectic form $\omega$. It produces another $DGA(k, \omega).$ Roughly speaking,
the origin of these data is the fact that de Rham cohomology of a smooth manifold with an
additional structure (complex, symplectic, homogeneous) carries a signature of this structure
upon its de Rham complex.

\smallskip

{\bf 3.6.1. Definition.} $(h,J)$ and  $(k,\omega )$ form a weak mirror pair iff these two DGAs are quasi--isomorphic.

\smallskip

{\bf 3.6.2. Proposition.} If $h$ and $k$ are nilpotent Lie algebras of common finite dimension,
then a homomorphism $DGA(h,J)\to DGA (k, \omega )$ is a quasi--isomorphism iff
it is an isomorphism.

\smallskip

In [Po], this is applied to the {\it extended deformations} of Kodaira surfaces
in the spirit of Merkulov: see [Me1], [Me2], and  our Sec. 5 below. Remarkably, it turns out
that in this world a Kodaira surface is its own mirror image.

\bigskip

\centerline{\bf 4. Extended deformation functors and controlling $L_{\infty}$--algebras}

\medskip

In the last two sections, we will sketch some extensions of  the controlling DGLAs
philosophy and constructions to the context of $\infty$--resolutions and higher categories.

\medskip

{\bf 4.1. $L_{\infty}$--algebras.} The notion of $L_{\infty}$--algebra, or homotopy Lie algebra $g$,
involves an infinite sequence of brackets on the $dg$--space $g$:
$$
\mu_n : \, \Lambda^n g \to g[2-n],\ n=1,2, \dots \infty\
\eqno(4.1)
$$
satisfying the relations, for all $n\ge 2$,
$$
\sum_{p+q =n+1} \sum_{\sigma\in {Sh^{-1}_{p,q}}} \roman{sgn} (\sigma )
 (-1)^{(p-1)q} (\mu_p\circ_1\mu_q)^{\sigma} =0 .
 \eqno(4.2)
 $$
We use here notations of  [LoVa], p.~365, Proposition 10.1.7, plus last line of the page, with $\mu_1=-d_g$. In particular, $Sh^{-1}_{p,q}$ denotes the set of
unshuffles, cf.~[LoVa], p.~16. 
 These conventions  agree also with those of [FiMaMar].
\medskip

{\bf 4.2. Maurer--Cartan equations for $L_{\infty}$--algebras.}  We put
for a $L_{\infty}$--algebra $g$:
$$
MC_{\infty} (g) := \{ x\in g^1 \,|\, \sum_{n=1}^{\infty} \frac{\mu_n (x^{\otimes n)}}{n!}=0 \}.
\eqno(4.3)
$$

{\bf 4.3. Homotopies in the set  $MC_{\infty} (g)$.}  The definition involving (2.2) can also be 
extended to this context, producing oriented paths between elements of $MC_{\infty}(g)$:
cf. [FiMaMar] and below.

\medskip

{\bf 4.4. Deligne  $\infty$--groupoids from $L_{\infty}$--algebras.} Generalising  sec.~2.5, consider an $L_{\infty}$--algebra
$g$ and 
  an Artin local algebra $A=\K\oplus m_A$.
  
  \smallskip
  
  Put 
 $$
  \roman{Ob}\, D_{\infty}(g, A) := MC_{\infty}(g\otimes m_A)
  \eqno(4.4)
  $$
  and for $x,y \in MC(g\otimes m_A)$
  $$
  \roman{Hom} (x,y) := \roman{paths\ from}\ x\ \roman{to}\ y.
  $$
  Actually, here we must not restrict ourselves by the composition of morphisms:
  equality between two compositions must be replaced by a path in the space of morphisms,
  and so on {\it ad infinitum}. So, as a functor of $A$, we will obtain
  an $\infty$--groupoid.
  \smallskip
  
  We omit here a formal description and instead treat a good motivating example from [FiMaMar].

\medskip

{\bf 4.5.  Semicosimplicial DGLAs.}  
Consider first the category $\Delta$ whose objects are finite sets
$$
[n] := \{0,1,\dots , n\},\ n=0,1,2,  \dots 
$$
and morphisms are order--preserving injective maps. Denote by
$$
\delta_{k,i} : [i-1] \to [i], \ k=0, \dots , i
$$
be the map with image $\{0,1,\dots , i\} \setminus \{k\}$. 

\smallskip

For a category $\Cal{X}$, call a {\it semicosimplicial $\Cal{X}$--object} any covariant
functor $\Delta \to \Cal{X}$. 

\smallskip

Thus a semicosimplicial DGLA  $g^{\Delta}$ is an  infinite sequence of DGLAs
$g_{i}$, $i=0,1, \dots \}$ connected by the morphisms $d_{k,i}: g_{i-1}\to g_i$
corresponding to $\delta_{k,i}$ and satisfying the same relations as  $\delta_{k,i}$.
\medskip

{\bf 4.6. Deligne groupoid $D(g^{\Delta}, A).$} The first step of its construction
leads to an infinite family consisting of objects $D(g_i,A)$ and respective morphisms, $A$ being
fixed.

\smallskip

The next step consists in passing to the homotopical limit. The groupoid
itself has as the objects ordered pairs of elements $\lambda,\mu \in (g_0^1\oplus g_1^0)\otimes m_A$ satisfying the conditions
$$
d\lambda +\frac{1}{2} [\lambda, \lambda]=0,\ e^{\mu}(d_{0,1}\lambda_0)=\lambda_1,
$$
$$
e^{d_{0,2}\mu}  e^{-d_{1,2}\mu} e^{d_{2,2}\mu}=1 .
$$
Finally, morphisms from  $\lambda_0,\mu_0$ to $\lambda_1,\mu_1$ 
are those elements $a\in g_0^0\otimes m_A$ for which
$$
e^a(\lambda_0)=\lambda_1,\  e^{-\mu_0} e^{-d_{1,1}a} e^{\mu_1} e^{d_{,1}a }=1.
$$

\medskip

{\bf 4.7.  A deformation problem.} In this subsection, $\K$ will denote an algebraically closed field of
characteristic zero. Consider a smooth algebraic variety $X$ over $\K$ and its
(finite) covering by open affine subsets $\U:= \{U_i\}$. It is well known that any infinitesimal
deformation over $A=\K\oplus m_A$ of an affine manifold $U$ is trivial, so they form a groupoid with
the single (isomorphism class of) object(s) $U\times Spec\, A$ and its automorphism
group $\roman{exp}\, (\Gamma (U, \T_U)\otimes_{\K}\, m_A)$. Thus (isomorphism classes of) all deformations of
$X$ over $A$ can be described as   the (noncommutative) cohomology set
$H^1( \U ,  \roman{exp}\, (\Gamma (U, \T_X)\otimes_{\K} m_A))$.

\smallskip

Extending this remark and building upon earlier work by E.~Getzler, V.~Hinich et al.
([Ge], [Hi1], [Hi2], [HiSch])
one can show ([FiMaMar]) that the whole \^Cech complex $C^*(\U,\T_X)$
has a natural structure of $L_{\infty}$--algebra and the whole $(\infty,1)$--groupoid
of deformations of $X$ is controlled by this $L_{\infty}$--algebra.

 \medskip
  
  {\bf 4.8. Further developments of the deformation theories and controlling DGLAs.}
  In Sec. 2 and 3,  I have explained the basics of controlling DGLAs philosophy as it was presented
  by the researchers of the 80s.  Bruno Vallette suggested me to include
  a brief picture of its development sketched in his message to me of July 22, 2017.
  With his permission, I reproduce below an edited version of his sketch.
  
  \smallskip
  
  The first remark concerns the initial Deligne formulation from [De]. The point is that later, when 
  $\infty$--groupoids were introduced, it became clear that the relevant DGLAs and 
  $L_{\infty}$--algebras are {\it filtered},
  and that $\infty$--groupoids are stable only wrt  filtered quasi--isomorphisms: see
  a modern treatment by Dolgushev--Rogers in arXiv:1407.6735 using   
   model category and homotopy arguments.

\smallskip
Furthermore, when Lurie (following the letter by Drinfeld of 1988 published as [Dr]) developed his version of the Deligne philosophy, he started with 
a generalisation of the notion of a general deformation functor. He then produced an infinity functor from DGLAs,
and an infinity functor in the opposite direction which together form
an $\infty$--equivalence.

\bigskip

\centerline{\bf 5. $F_{\infty}$--structures on extended deformation spaces}

\medskip

{\bf 5.1. $F$--manifolds.} We start with description of a class of 
manifolds whose tangent sheaf is endowed with (commutative, associative) multiplication
such as (1.2) in the section 1 above.

\smallskip

Below $M$ denotes a (super)manifold: it can be $C^{\infty}$, or $An$, or (partly) formal,
eventually with odd (anticommuting) coordinates.
The ground field is denoted
$K$, usually we choose $K=\bold{C}$.

\smallskip

The structure sheaf is denoted $\Cal{O}_M$, the tangent sheaf $\Cal{T}_M$.
The tangent sheaf is a locally free $\Cal{O}_M$--module;
its (super)rank is called the (super)dimension of $M$.

\smallskip

Now start with  a linear $K$--(super)space $A$ endowed with $K$--bilinear
commutative multiplication and a $K$--bilinear Lie bracket.

\smallskip

The Poisson tensor of such a structure  assigns to $a,b,c \in A$ the element
$$
P_a(b,c):=[a,bc]-[a,b]c-(-1)^{ab}b[a,c].
$$
This definition can be easily extended to sheaves.

\smallskip

For a manifold $M$ as above, $\Cal{O}_M$ has a natural commutative multiplication, whereas  $\Cal{T}_M$
has a natural Lie structure. 
\smallskip
Poisson structure involves introducing additional Lie structure upon  $\Cal{O}_M$, whereas
$F$--structure 
involves  introducing additional multiplication  upon $\Cal{T}_M$, satisfying axioms below.
Below we compare the axioms and particular cases of these two structures.

\newpage

\centerline{}

\bigskip

{\bf POISSON STRUCTURE}        \quad\quad\quad\quad\quad\quad\quad                   {\bf F--STRUCTURE}

 \medskip

 $K$--bilinear (super)Lie    \quad\quad\quad\quad\quad\quad\quad\quad\  \                       $\Cal{O}_M$--bilinear associative, commutative

  bracket $\{,\}$ on $\Cal{O}_M$ \quad\quad\quad\quad\quad\quad\quad\quad\quad \,             unital multiplication $\circ$ on $\Cal{T}_M$

        satisfying identity           \quad\quad\quad\quad\quad\quad\quad\quad\quad\quad                      satisfying ``$F$-- identity'':

     \smallskip

 $ P_f(g,h) \equiv 0$   \quad\quad\quad \quad\quad\quad\quad\quad\quad\quad\quad\quad\quad\  \,       $P_{X\circ Y}=X\circ P_Y+(-1)^{XY}Y\circ P_X$

  \bigskip  

 Equivalently: each local    \quad\quad\quad\quad\quad\quad\quad\ \                              Each local vector field on $M$

  function $f$ on $M$ becomes a local  \quad\quad\quad\quad  becomes a local function  
  
  vector field $X_f$ on $M$: \quad\quad\quad\quad\quad\quad\quad\quad\quad  on the spectral cover

   $X_f(g):=\{f,g\}$  \quad\quad\quad\quad\quad\quad\quad\quad\quad\quad\quad\quad\                        $\widetilde{M}:=Spec_{\Cal{O}_M}(\Cal{T}_M,\circ )$

 \medskip

 {\it Special case (local):}   \quad\quad\quad\quad\quad\quad\quad\quad\quad\quad\quad                           {\it Special  case (local):}

   \smallskip                                                                                                              

 Symplectic structure    \quad\quad\quad\quad\quad\quad\quad\quad\quad\                                                                 Semisimple $F$--structure

  in canonical coordinates:     \quad\quad\quad\quad\quad\quad\quad\ \                                in Dubrovin's canonical 

 $\{f,g\}=\sum_{i=1}^n\partial_{q_i}f\partial_{p_i}g  -\partial_{q_i}g\partial_{p_i}f$ \quad\quad\quad\quad\quad\ \   coordinates $u^a$:
                                
   No local paramters, but       \quad\quad\quad\quad\quad\quad\quad\quad                               $\partial_a\circ \partial_a =\delta_{ab}\partial_a$

  a large symplectomorphism group.    \quad\quad\quad\                                                 No local  parameters;

 \quad\quad\quad \quad\quad\quad \quad\quad\quad\quad\quad\quad \quad\quad\quad \quad\quad\quad\quad local automorphisms 

 \quad\quad\quad \quad\quad\quad \quad\quad\quad\quad\quad\quad \quad\quad\quad \quad\quad\quad\quad $u^a\mapsto u^{\sigma (a)}+c^a$

\bigskip

In a recent article [Do], there is very interesting description of the operad $\roman{\bold{FMan}}$, classifying
algebras $(A, \circ , [,])$ whose basic operation $\circ$ is commutative and associative,
basic operation $[,]$ is the Lie bracket, and finally, their compatibility is expressed by
the $F$--identity. Notice that $F$--identity is a cubic one in the operadic sense
so that the connection between $\roman{\bold{FMan}}$ and quadratic operads
Associativity, Poisson and Pre--Lie operads is quite surprising.

\smallskip

Since among these three operads the last one is less well known, we briefly
recall that the Pre--Lie operad classifies pre--Lie  algebras, and the latter
are defined by binary product whose associator is right symmetric: see [LoVa],
Sec~13.4. In [DoShVa2], a version of Deligne groupoid and pre--Lie deformation formalism is 
developed.

\medskip

{\bf 5.2. Geometric meaning of the $F$--identity} ([HeMaTe]).
For any (super)manifold $M$, consider
the sheaf of those functions on the 
cotangent manifold $T^*M$  which are polynomial
along the fibres of projection $T^*M\to M$: 
that is, the relative symmetric algebra $Symm_{\Cal{O}_M}(\Cal{T}_M)$.

\medskip

It is a sheaf of $\Cal{O}_M$--algebras, multiplication in which we denote $\bold{\cdot}$

\smallskip
Consider now a triple $(M,\circ, e)$ where $\circ$ is a commutative associative $\Cal{O}_M$--bilinear  multiplication
on $\Cal{T}_M$, eventually with identity $e$.

\medskip

There is an obvious homomorphism of $\Cal{O}_M$--algebras 
$$
(Symm_{\Cal{O}_M}(\Cal{T}_M), \cdot )\to (\Cal{T}_M, \circ )
$$

{\bf 5.2.1. Theorem}.  {\it The multiplication $\circ$ satisfies the $F$--identity
$$
P_{X\circ Y}=X\circ P_Y+(-1)^{XY}Y\circ P_X
$$
iff its kernel is stable with respect to 
the canonical Poisson brackets on $T^*M$.}
\medskip
 In other words, $F$--identity is equivalent to the fact that the spectral cover of $M$
considered as a closed subspace of its cotangent bundle is coisotropic of
maximal 
dimension.

\medskip

 NB. The spectral cover $\widetilde{M}:=Spec_{\Cal{O}_M}(\Cal{T}_M,\circ )$ of $M$ is not necessarily
a manifold.
Its structure 
sheaf may have zero divisors and nilpotents.

\smallskip

However, it is a manifold, if the $F$--manifold $M$ is semisimple.
\smallskip
 Conversely, an embedded submanifold $N\subset T^*M$ is
the spectral cover of some 
semisimple 
$F$--structure iff $N$ is Lagrangian.

\medskip

{\bf 5.3. Local decomposition theorem.}
 Sum of two $F$--manifolds is defined by:
$$
(M_1,\circ_1,e_1) \oplus (M_2,\circ_2,e_2) :=
(M_1\times M_2,\circ_1\boxplus\circ_2,e_1\boxplus e_2)
$$

A manifold is called {\it indecomposable} if it cannot be represented
as a sum in a nontrivial way.

\smallskip
For any point
$x$ of a pure even $F$--manifold $M$, the tangent space $T_xM$
is endowed with the 
structure of a commutative finite dimensional $K$--algebra.
This $K$--algebra can be represented  
as the direct sum of local
$K$--algebras. The decomposition is unique in the following sense:
the set of 
pairwise orthogonal idempotent tangent
vectors determining is well defined.

\bigskip

{\bf 5.3.1. Decomposition Theorem.}
 {\it Every germ $(M,x)$ of a
complex analytic $F$--manifold  
decomposes 
into a direct sum
of indecomposable germs such that for each summand,
the tangent algebra 
at $x$ is a local algebra.

\smallskip

This decomposition is unique
in the following sense: the set of pairwise orthogonal idempotent
vector fields determining it is well defined.}
\medskip

{\bf 5.3.2. Comments.}
 (i) If $(T_xM,\circ )$ is semisimple, this
theorem is equivalent to the existence 
(and 
uniqueness) of Dubrovin's coordinates.

\smallskip

(ii) A proof of this theorem is based upon interpretation of the
basic identity of the $F$--structure as integrability 
condition.

\smallskip

(iii) For $F$--manifolds with a compatible flat structure,
there exists a considerably more 
sophisticated operation 
of {\it tensor product} which we omit here.

\smallskip

Furthermore, we have ([He], Theorems 5.3 and 5.6):

\medskip

{\bf 5.4. Theorem.}
  {\it (i) The spectral cover space $\widetilde{M}$
of the canonical $F$--structure on the germ of the 
unfolding space
of an isolated hypersurface singularity is smooth.

\smallskip

(ii) Conversely, let $M$ be an irreducible germ of a generically
semisimple $F$--manifold with 
the smooth spectral cover $\widetilde{M}.$
Then it is (isomorphic to) the germ of the unfolding space
of an 
isolated hypersurface singularity. Moreover, any isomorphism
of germs of such unfolding 
spaces compatible with their
$F$--structure comes from a {\it stable right equivalence} of the
germs 
of the respective singularities.}

\medskip

Recall that the stable right equivalence is generated by adding sums
of squares of coordinates 
and making invertible local analytic
coordinate changes.

\smallskip

In view of this result, it would be important to understand
the following

\medskip

{\bf 5.4.1. Problem.} {\it Characterize those varieties $V$ for which the genus zero
 quantum cohomology Frobenius 
 spaces $H^*_{quant}(V)$ have
smooth spectral covers.}

\medskip

Theorem 5.4 above produces for such manifolds a weak version
of Landau--Ginzburg model, and 
thus gives a partial
solution of the mirror problem for them.

\medskip

{\bf 5.5. From $F$--manifolds to Frobenius manifolds.} We start with an incomplete 
description of such a passage and steps involved in it.

\smallskip
 A Frobenius manifold is an
$F$--manifold endowed 
with a compatible flat structure $\nabla$, an Euler vector  field $E$
and a (pseudo)--Riemannian 
metric $g:\,S^2(\Cal{T}_M)\to \Cal{O}_M$
such that

\medskip

(i) $g$ is flat, and $\nabla$ = the Levi--Civita connection of $g$.

\smallskip

(ii) $g(X\circ Y,Z)= g(X, Y\circ Z).$

\smallskip 

An Euler field $E$ is compatible with Frobenius structure if

\smallskip

(iii) $Lie_Eg=Dg$ for a constant $D$.

\medskip

NB. This is only an incomplete version because not all restrictions
and 
compatibility conditions on extra structures  are spelled out explicitly below.
\medskip
One condition of compatibility is stated below in more detail because it will be important 
for the definition of $F_{\infty}$--structure.

\medskip

{\bf 5.6. Compatible flat structures.}
  An (affine) flat structure on a (super) manifold  $M$ is given by any
of the following equivalent data:

\medskip

(i) A torsionless flat connection $\nabla_0:\,\Cal{T}_M\to
\Omega^1_M\otimes_{\Cal{O}_M}\Cal{T}_M $.

\smallskip

(ii) A local system $\Cal{T}_M^f\subset \Cal{T}_M$ of flat vector fields,
which forms a sheaf of (super)commutative 
Lie algebras of rank $\roman{dim}\,M$
such that $\Cal{T}_M=\Cal{O}_M\otimes \Cal{T}_M^f$.

\smallskip

(iii) An atlas whose transition functions are affine linear.

\bigskip

Assume that $\Cal{T}_M$ is endowed with an $\Cal{O}_M$--bilinear
(super)commutative and associative 
multiplication $\circ$,
and eventually with unit $e$. 

\medskip

NB $F$--identity is not yet postulated! 

\medskip

{\bf 5.6.1. Definition} ([Man2]).  {\it a) A flat structure $\Cal{T}_M^f$
on $M$ is called compatible with $\circ$, 
if in a neighborhood of any point there exists a vector field
$C$ such that for arbitrary 
local flat vector fields $X,Y$
we have
$$
X\circ Y= [X,[Y,C]].
$$
$C$ is called then a local vector potential for $\circ .$

\smallskip
b)  $\Cal{T}_M^f$ is called compatible with $(\circ ,e)$, if
a) holds and moreover,  $e$ is flat.
 
\medskip

{\bf 5.6.2. Proposition.}
 {\it If $\circ$ admits a compatible flat structure,
then it automatically 
satisfies the $F$--identity.}}

\medskip

{\bf 5.7. $F_{\infty}$--manifolds.}  The $\infty$--version of $F$--manifolds discussed below
generalises only  the case of infinitesimally deformed germ of a manifold $(M,*)$ and replaces it
by its formal {\it smooth} graded dg resolution 
$(\M , *)$  supplied with a smooth degree one vector field $\partial$ satisfying
$$
[\partial ,\partial ] =0,\quad \partial I \subset I^2
$$
where $I$ is the ideal of $*$.

\smallskip

The role of $\circ$ on $(\M, *,\partial )$
will now be played by
a structure of $C_{\infty}$--algebra
$$
\mu_{\bullet} = \{\mu_n\}_{n\ge 1} : \otimes_{\O_\M} \T_{O_\M} \to \T_{O_\M}
$$

\smallskip

The former $F$--identity in this context becomes the first step
of the ladder:  $[\mu_2,\mu_2]=0$ where by definition $[\mu_2,\mu_2]:\ \otimes^4_{\O_M}\T_M \to \T_M$ is given by 
$$
[\mu_2,\mu_2] (X,Y,Z,W) := the\ left\ hand\ side\ of\ F-identity
$$

The whole ladder involves  a system of 
``polybrackets'' $[\mu_{\bullet}, \mu_{\bullet}]^{\nabla}$ depending on the additional choice
of a torsion free affine connection $\nabla$
and subsequent passage to the its cohomology class
$$
[[\mu_{\bullet}, \mu_{\bullet}]] \in H( \otimes_{\O_\M}^{\bullet}\T_{\M}^*\otimes_{\O_\M} \T_{\M})
$$

The $F_\infty$--identity then reads
$$
[[\mu_{\bullet}, \mu_{\bullet}]]=0,
$$
and it defines a structure of $F_{\infty}$--manifold upon $(\M,*,\partial )$.

\medskip

{\bf 5.8. Theorem.} ([Me2]).
 {\it (i) The formal dg manifold associated with
the Hochschild 
cohomology of an associative algebra is an $F_{\infty}$--manifold.

\smallskip

(ii) The formal dg manifold associated with singular cohomology of  a compact 
topological
space is an $F_{\infty}$--manifold.}

\medskip

In [DoShVa1], these results are somewhat generalised and/or strengthened.

\bigskip

{\bf Acknowledgements.} This survey was written as a contribution for the collection 
``Handbook of Mirror Symmetry for Calabi--Yau manifolds
and Fano manifolds'' (in preparation), by invitation of  Shin--Tung Yau.
Its conception was greatly influenced by  remarks by Jean--Pierre Serre that
he likes in mathematical papers ``precision combined with informality'' and ``side remarks, open problems,
and such'' that very often are ``more interesting than the theorems actually proved''
(interview in Singapore, 1985).
The initial version of the text was written as a presentation of my talk at
the MPIM, Bonn Arbeitstagung 2017 ``Physical Mathematics''  organised by Peter Teichner.
Subsequent drafts were very carefully read by Bruno Vallette, who made many
useful suggestions incorporated in the text, and by Sergey Merkulov whose work
furnished a very strong motivation for the paper.
\smallskip

I am deeply grateful to all these people.

\bigskip
\centerline{\bf References}

\medskip

[Ar] Arnold, V. I. {\it Normal forms of functions near degenerate critical points, the Weyl groups $A_k,D_,E_k$
and Lagrangian singularities (Russian)}. Funkc. Anal. i Prilozhen.  vol.~6, no.~4, 1972, pp.~ 3--25. 

\smallskip

[Art1] M.~Artin. {\it On solutions of analytic equations.} Inv.~Math. 5, 1968, pp.~277--291.

\smallskip

[Art2] M.~Artin. {\it Lectures on Deformations of Singularities.} Lectures on Math. and Phys. 
Tata Institute 54 (1976).

\smallskip


\smallskip
[BaMan] A.~Bayer, Yu.~Manin. {\it (Semi)simple exercises in quantum cohomology.} In: The
Fano Conference Proceedings, 
ed. by A. Collino, A. Conte, M. Marchisio, Universit\`a
    di Torino, 2004, pp.~143--173. arXiv: math.AG/0103164

\smallskip

[BehMan] K.~Behrend, Yu.~Manin. {\it Stacks of stable maps and Gromov--Witten invariants.}
Duke Math. Journ., 85:1, 1996, pp.~1--60.
 
\smallskip
[BelGeKo] A.~Belavin, D.~Gepner, Y.~Kononov. {\it Flat coordinates for Saito
Frobenius manifolds and String theory.} arXiv:1510.06970

\smallskip

[BuMi] R.--O.~Buchweitz, J.~Millson. {\it CR--Geometry and Deformations of Isolated
Singularities.} Memoirs of AMS, Vol.~125, N.~597, 1997, 96 pp.

\smallskip

[ChMPe] P.~E.~Chaput, L.~Manivel, N.~Perrin. {\it Quantum Cohomology of minuscule
homogeneous spaces.} Transformation Groups, vol.~13, no.~1, 2008, pp.~47--89.

\smallskip

[ClLaPo] R.~Cleyton, J.~Lauret, Yat Sun Poon. {\it Weak Mirror Symmetry of Lie
Algebras.}  J.~Symplectic Geom.  vol.~8, no. 1, 2010, pp.~37--55,             arXiv:1004.3264. 22pp.

\smallskip

[ClPo] {\it Differential Gerstenhaber Algebras Associated to Nilpotent Algebras.} Asian J.~Math. vol.~12, no. 2,
2008, pp.~255--249.

\smallskip

[ClOvPo11]  R.~Cleyton, G.~Ovando, Yat Sun Poon.  {\it  Weak Mirror Symmetry of Complex 
Symplectic Lie Algebras.} J.~Geom.~Phys. 61(2011), no. 8, 2001, pp.~1553--1563.     arXiv:0804.4787. 

\smallskip

[De] P.~Deligne. {\it Letter to J.~J.~Millson.} April 24, 1986.

\smallskip

[Do]  V.~Dotsenko. {\it Algebraic structures of $F$--manifolds via
pre--Lie algebras.} arXiv:1706.07340

\smallskip


\smallskip

[DoShVa1] V.~Dotsenko, S.~Shadrin, B.~Vallette. {\it De Rham cohomology and homotopy Frobenius
manifolds.} arXiv:1203.5077

\smallskip

[DoShVa2] V.~Dotsenko, S.~Shadrin, B.~Vallette. {\it Pre--Lie deformation theory.}
arXiv:1502.03280. 31 pp.

\smallskip

[Dr] V.~Drinfeld. {\it A letter from Kharkov to Moscow.} EMS Surv.~Math.~Sci.,
2014, vol.~1, no.~2, pp.~241--248.

\smallskip

[FiMaMar] D.~Fiorenza, M.~Manetti, E.~Martinengo. {\it Cosimplicial DGLAs in Deformation Theory.}
Comm. Algebra 40, 2012, vol. 6, pp.~2243--2260

arXiv:math0803.0399

\smallskip

[Ge] E.~Getzler. {\it Lie theory for nilpotent $L_{\infty}$--algebras.} Annals of Math.,
170 (1), 2009, pp.~271--301. arXiv:math/0404003.

\smallskip

[GoMi] W.~M.~Goldman, J.~J.~Millson. {\it The deformation theory of
representations of fundamental groups of compact K\"ahler manifolds.}
in: Publ.~Math. IHES, tome 67, no.~2, 1988, pp.~43--69.

\smallskip

[GoPe] V.~Gorbunov, V.~Petrov. {\it Schubert calculus and singularity theory.} 

arXiv:1006.1464

\smallskip

[Gr] H.~Grauert. {\it \"Uber die Deformationen isolierter Singularit\"aten analytischer Mengen.}
Inv.~Math. 15(1972), 171--198.

\smallskip

[He] C.~Hertling. {\it Frobenius spaces and and moduli spaces for singularities.}
Cambridge University Press, 2002, ix+270 pp.

\smallskip

[HeMan] C.~Hertling, Yu.~Manin. {\it Unfoldings of meromorphic connections and a construction of Frobenius manifolds.} In: Frobenius manifolds, 113--144, Aspects Math., E36, Friedr. Vieweg, Wiesbaden, 2004. 


\smallskip

[Hi1] V.~Hinich. {\it Descent of Deligne groupoids.}  IMRN 1997, no.~5, pp.`223--239
arXiv:alh--geom/9606010 . 

\smallskip

[Hi2]  V.~Hinich. {\it Deformations of sheaves of algebras.}, Advances in Math., vol.~195, 2005,
pp.~102--164.

\smallskip

[HiSch] V.~Hinich, V.~Schechtman. {\it Deformation Theory and Lie Algebra Homology I,II.}
Algebra Colloquium 4:2 (1997), pp. 213--240 and 4:3 (1997), pp. 291--316.

\smallskip

[KoMan]  M.~Kontsevich, Yu.~Manin. {\it Gromov--Witten classes, quantum cohomology, and enumerative geometry.}
Comm. Math. Phys., 164:3 (1994), pp.~525--562
\smallskip

[LamTem] Th.~Lam, N.~Templier. {\it The Mirror Conjecture for minuscule
flag varieties.} arXiv:1705.00758.

\smallskip

[LoVa] J.~L.~Loday, B.~Vallette. {\it Algebraic Operads.} Springer, 2012, xxiv + 634 pp.

\smallskip

[Lu] J.~Lurie. {\it Derived Algebraic Geometry X: Formal Moduli Problems.} 

www.math.harvard.edu/~lurie
\smallskip
[Ma1] M.~Manetti. {\it Deformation theory via differential graded Lie algebras.} arXiv:math/0507284
\smallskip

[Ma2] M.~Manetti. {\it Extended deformation functors.} IMRN no. 14, 2002, pp. 719--756.
arXiv: 9910071
\smallskip
[Ma3] M.~Manetti. {\it On some formality criteria for DG--Lie algebras.} J.~Algebra 438, 2015,
pp. 90--118. arXiv:1310.3048

\smallskip

[Man1]  Yu.~Manin.  {\it Frobenius manifolds, quantum cohomology, and moduli
spaces.}  AMS Colloquium Publications, vol. 47, Providence, RI, 1999,
xiii+303 pp.

\smallskip

[Man2] Yu.~Manin. {\it Manifolds with multiplication on the tangent sheaf.}
 Rendiconti Mat. Appl., Serie VII, vol.26 (2006), 69--85.  arXiv:0502578

\smallskip

[ManSm1] Yu.~Manin, M.~Smirnov. {\it On the derived category of $\overline{M}_{0,n}$}.
 Izvestiya of Russian Ac. Sci., vol. 77, No 3, 2013, 93--108.
 arXiv:1201.0265

\smallskip

[ManSm2] Yu.~Manin, M.~Smirnov. {\it Towards motivic quantum cohomology of $\overline{M}_{0,S}$.}
Proc. of the Edinburg Math. Soc., Vol. 57 (ser. II), no 1, 2014,
pp. 201--230. Preprint arXiv:1107.4915

\smallskip

[MatY] J.~N.Mather, S.~S.~T.~Yau.  {\it Classification of isolated hypersurface singularities 
by their moduli algebras.}  Inv.~Math., vol.~69, 1982, pp.~243--251.

\smallskip

[Me1] S.~Merkulov. {\it $Frobenius_{\infty}$  invariants of homotopy Gerstenhaber algebras I.}
Duke Math. J. 105, 2000, pp. 411--461.

\smallskip

[Me2] S.~Merkulov. {\it Operads, deformation theory and $F$--manifolds.} In: Frobenius Manifolds
(ed.~by C.~Hertling, M.~Marcolli).  Vieweg, 2004, pp.~213--251.
arXiv:02100478.
\smallskip

[Mi] J.~Millson. {\it Rational homotopy theory and deformation problems from algebraic
geometry.} Proc.~ICM 90, Kyoto, Math.~Soc.~Japan, 1991, pp.~549--558.

\smallskip

[Po] Yat Sun Poon. {\it Extended Deformation of Kodaira Surfaces.} 
J.~reine angew.~Math., 590 (2006), pp. 45--65.
arXiv:math/0402440

\smallskip

[Pr] J.~P.~Pridham. {\it Unifying derived deformation theories.} Advances in Math., 224 (3),
2010, pp.~772--826.

{\it Corrigendum: } Advances in Math., 228, 2011, pp. 2554--2556.

\smallskip

[R-N] D.~Robert--Nicoud. {\it Representing the Deligne--Hinich--Getzler $\infty$--groupoid.}
arXiv:1702.02529, 13 pp.

\smallskip

[Sa] K.~Saito. {\it Primitive forms for universal unfolding of a function with an isolated critical
point.} Journ.~Fac.~Sci.~Univ.~Tokyo, sec IA, vol.~28, 1981, pp.~775--792.

\smallskip

[Sc] J.~Scherk. {\it A propos d'un th\'eor\`eme de Mather and Yau.} C.~R.~Ac.~Sci.~Paris, S\'er.~I, t.~296,
1983, pp.~513--515.

\smallskip

[Sch] M.~Schlessinger. {\it Functors of Artin rings.} Trans.~AMS, vol. 130, 1968, pp.~ 208--222.

\smallskip 

[SchSt]  M.~Schlessinger, J.~Stasheff. {\it The Lie algebra structure
of tangent cohomology and deformation theory.} J.~Pure Appl. Algebra, vol. 38, no. 2--3, 1985, pp. 313--322.

\smallskip

[To] B.~To\"en. {\it Probl\`emes de modules formels [d'apr\`es V.~Drinfeld, V.~Hinich, M.~Kontsevich,
J.~Lurie \dots].} S\'em. Bourbaki, Jan.~16, 2016.
 https://perso.math.univ-toulouse.fr/btoen/files/2012/04/Bourbaki-Toen-2016-final1.pdf .

\bigskip

\enddocument